\documentclass[12pt]{article}
\usepackage{amsmath, amsfonts, epsfig, txfonts}
\def\theequation{\arabic{equation}}

   \def\theequation{\thesection.\arabic{equation}}
   \def\theequation{\arabic{equation}}

\newsavebox{\savepar}

\newtheorem{remark}{Remark}

\newtheorem{example}{Example}
\newtheorem{definition}{Definition}
\newtheorem{theorem}{Theorem}

\def\appendix{\par
\setcounter{section}{0}
   \setcounter{equation}{0}
    \def\@chapapp{APPENDIX}
    \def\thesection{\Alph{section}}
    \def\theequation{A.\arabic{equation}}
    \def\section{
    \refstepcounter{section}
    \@startsection {section}{A}{0pt}{-3.5ex \@plus -1ex \@minus -.2ex}
                   {0.3ex \@plus.2ex}{\normalsize\sffamily\bfseries} }}

\def\eop{\hbox{{\vrule height7pt width3pt depth0pt}}}

\makeatletter
\newcommand{\least}{\let\CS=\@currsize\renewcommand{\baselinestretch}{1}\tiny\CS}

\newcommand{\oneandahalfspacing}{\let\CS=\@currsize\renewcommand{\baselinestretch}{1.2}\tiny\CS}
\newcommand{\doublespacing}{\let\CS=\@currsize\renewcommand{\baselinestretch}{2.5}\tiny\CS}
\usepackage{pst-plot}
\usepackage{amsmath}
\usepackage{amsfonts}
\usepackage{euscript}
\usepackage{oldgerm}

\usepackage{rotating}

   \renewcommand{\baselinestretch}{1.3}

\begin{document}

\newcommand{\namelistlabel}[1]{\mbox{#1}\hfil}
\newenvironment{namelist}[1]{%
\begin{list}{}
{
\let\makelabel\namelistlabel
\settowidth{\labelwidth}{#1}
\setlength{\leftmargin}{1.1\labelwidth}
}
}{%
\end{list}}

\newcommand{\be}{\begin{equation}}
\newcommand{\ee}{\end{equation}}
\newcommand{\dist}{{\rm\,dist}}
\newcommand{\sspan}{{\rm\,span}}
\newcommand{\re}{{\rm Re\,}}
\newcommand{\im}{{\rm Im\,}}
\newcommand{\sgn}{{\rm sgn\,}}
\newcommand{\beano}{\begin{eqnarray*}}
\newcommand{\eeano}{\end{eqnarray*}}
\newcommand{\bea}{\begin{eqnarray}}
\newcommand{\eea}{\end{eqnarray}}

\newcommand{\ba}{\begin{array}}
\newcommand{\ea}{\end{array}}
\newcommand{\hone}{\mbox{\hspace{1em}}}
\newcommand{\hon}{\mbox{\hspace{1em}}}
\newcommand{\htwo}{\mbox{\hspace{2em}}}
\newcommand{\hthree}{\mbox{\hspace{3em}}}
\newcommand{\hfour}{\mbox{\hspace{4em}}}
\newcommand{\von}{\vskip 1ex}
\newcommand{\vone}{\vskip 2ex}
\newcommand{\vtwo}{\vskip 4ex}
\newcommand{\vthree}{\vskip 6ex}
\newcommand{\vfour}{\vspace*{8ex}}
\newcommand{\norm}{\|\;\;\|}
\newcommand{\integ}[4]{\int_{#1}^{#2}\,{#3}\,d{#4}}
\newcommand{\inp}[2]{\langle {#1} ,\,{#2} \rangle}
\newcommand{\vspan}[1]{{{\rm\,span}\{ #1 \}}}
\newcommand{\R} {{\mathbb{R}}}

\newcommand{\B} {{\mathbb{B}}}
\newcommand{\C} {{\mathbb{C}}}
\newcommand{\N} {{\mathbb{N}}}
\newcommand{\Q} {{\mathbb{Q}}}
\newcommand{\LL} {{\mathbb{L}}}
\newcommand{\Z} {{\mathbb{Z}}}

\newcommand{\BB} {{\mathcal{B}}}
\newcommand{\dm}[1]{ {\displaystyle{#1} } }
\def \stackt{{\stackrel{.}{.\;.}\;\;}}
\def \stackb{{\stackrel{.\;.}{.}\;\;}}
\def \olu{\overline{u}}
\def \olv{\overline{v}}
\def \olx{\overline{x}}
\def \olp{\overline{\partial}}
\def\diag{{\;{\rm diag } \; }}
\thispagestyle{empty}

\newcommand{\alert}[1]{\fbox{#1}}

\begin{center}
{\large \bf  Average    Collapsibility of Some  Association Measures}\\
\end{center}

 \vtwo

\noindent {\bf {P. Vellaisamy}}\\
{\it Department of Mathematics, Indian Institute of Technology
Bombay, Mumbai-400076, India.}

\vtwo

\noindent {\bf Abstract.} Collapsibility deals with the conditions
under which a conditional (on a covariate $W$) measure of
association between two random variables $X$ and $Y$ equals the
marginal measure of association, under the assumption of
homogeneity over the covariate. In this paper, we discuss the
average collapsibility of certain well-known measures of
association, and also with respect to a new measure of association.
The concept of average   collapsibility is more general than
collapsibility, and requires that the conditional average
 of an association measure equals the corresponding marginal measure. Sufficient conditions for the
  average   collapsibility of the measures
under consideration are obtained. Some difficult,  but interesting, counter-examples
are constructed. Applications to linear, Poisson,
 logistic and negative binomial regression models are addressed. An extension to the case of multivariate covariate $W$ is also discussed.

\vtwo

\vone {\noindent \it Key words:} { Average  collapsibility,
collapsibility, conditional distributions, linear and non-linear
regression models, measures of association, Yule-Simpson paradox.}

\vone {\noindent \bf {AMS Subject Classification}:} Primary:
62H20; ~Secondary: 62J05, 62J12.

\vone
\noindent {\bf 1. Introduction}

\noindent The study of association between two random variables arises in several applications. Several measures, nonparametric
 in nature, have been proposed in the literature. Often, the random variables of interest, say $X$ and $Y$,
 may be associated because of their association
 with another variable $W$, called a covariate or a background variable.
In this case, we need to investigate the conditional association
measure between $X$ and $Y$ given $W$, and compare it with  the
marginal association measure between $X$ and $Y$. It is in general
possible that the conditional association measure
 may be positive, while the marginal association measure may be negative. Such an effect reversal is called the
 Yule-Simpson paradox attributed to Yule (1903) and Simpson
 (1951). When Yule-Simpson paradox or the effect reversal does not occur, and a
 conditional measure of association equals the marginal measure,
 we say that the measure is collapsible over the covariate $W$. Collapsibility is an important issue associated with  data analysis, analysis of contingency tables, causal inference, regression analysis, epidemiological studies and the design of experiments; see, for example, Cox and
 Wermuth (2003), Ma \emph{ et al.~}(2008), and Xie \emph{ et al.~}(2008) for applications and discussions. \\

 There have been several notions of collapsibility,
namely, simple, strong and uniform collapsibility. These issues
have been addressed in several different contexts such as the
analysis of contingency tables, regression models and association
measures; see for example Bishop (1971), Cox (2003), Cox and
Wermuth (2003), Geng (1992), Ma {\em et~al.} (2006), Vellaisamy
and Vijay (2008), Wermuth (1987, 1989),  Whittemore (1978), and Xie
{\em{et al.~}}(2008). Cox and Wermuth (2003) studied the concept
of distribution dependence and discussed the conditions under
 which no effect reversal occurs. Xie \emph{ et al.~}(2008) discussed the simple collapsibility and the uniform
 collapsibility of the following association measures :
\begin{eqnarray}
 \begin{array}{lll}
 (i) & \displaystyle \frac{\partial}{\partial x}E(Y\mid x) & {\rm (expectation ~ dependence)}\nonumber \\
 (ii) & \displaystyle \frac{\partial^2}{\partial x\partial y}\log f(x,y) & {\rm (mixed ~ derivative ~ of ~
 interaction)}\nonumber \\
 (iii) & \displaystyle \frac{\partial}{\partial x}F(y\mid x) & {\rm (distribution ~ dependence)}.\end{array}\nonumber
\end{eqnarray}
They discussed also the stringency of the above measures for positive
association, studied the conditions for no effect
 reversal (after marginalization over $W$) and  obtained  the necessary and sufficient conditions for uniform
 collapsibility
 of mixed derivative of interaction, among other results.
Recently, Vellaisamy (2011) introduced a new concept of average
collapsibility and discussed it with respect to the distribution
dependence and  the quantile regression
 coefficients. It is shown that average   collapsibility is a general concept and coincides with  collapsibility
 under the condition of homogeneity. In the same spirit, we discuss in this paper the average   collapsibility of
  expectation dependence, and mixed
 derivative of interaction measures which have relevance to linear and logistic regression models.
  Also, a new measure of
 association, namely,
\begin{eqnarray*}
 ~~~~~(iv) ~~ \displaystyle \frac{\partial}{\partial x} \log E(Y\mid x) ~~~~  {\rm (\log expectation ~ dependence)}
\end{eqnarray*}
 is introduced and
 its average   collapsibility conditions are investigated. This measure has a direct application to
 Poisson and negative binomial regression  models. In the  last
 section, some results are extended to the case of multivariate
 covariate $W.$

\vone
 \noindent {\bf 2. The Average  Collapsibility Results}

\noindent Let $(Y, X,W)$ be a random vector, where our interest is mainly on the association between $Y$ and $X$, and $W$ is treated as a covarite. We assume for simplicity that $X$ and $W$ are continuous, unless stated otherwise. Note that $Y$ has a monotone (increasing) regression function of $X$ if $E(Y|X=x)$ is increasing in $x$ or equivalently the expectation dependence function (EDF) ${\partial}E(Y\mid x)/{\partial x} \geq 0$. We first discuss the average collapsibility results for the $EDF$ and introduce the
following definition.

 {\definition
The expectation dependence function ($EDF$) is average collapsible
over $W$ if
\begin{equation} \label{eqn2.1}
E_{W|x} \left(\frac{\partial}{\partial x}E (Y|x,W)\right) =
\frac{\partial}{\partial x}E(Y|x),\hspace{4mm}   \mbox{for all} ~
~x .
\end{equation}}

 The following result gives sufficient conditions for the average
collapsibility  of $EDF.$ In the sequel,  $X \Perp Y$ and $X \Perp
Y |W$ respectively denote the independence of $ X$ and $Y$, and
the conditional independence  of $ X$ and $Y$ given $W$. We assume henceforth
all the partial derivatives exist and are continuous so that that the differentiation and
integration can be interchanged.
\begin{theorem} \label{thm2.1}  The $EDF$    $\dfrac{\partial}{\partial x}E(Y|x,w)$ is average  collapsible
over $W$  if either \\
$ \hspace*{3cm}(i)~ E(Y|x,w)$ is independent of $w$, or \\
$\hspace*{3cm}(ii)~ X\Perp W$ \\
holds.
\end{theorem}

The condition that $E(Y|x, w)$ is independent of $w$ implies the
homogeneity of $EDF$ and in this case both uniform collapsibility
(Part (a) of Theorem 3.4 of Xie \emph{et al.~}(2008)) and average
collapsibility hold.  However, when the $EDF$ is not homogeneous over $w$,
average collapsibility may still hold if (and only if) $ X \Perp W$.
 Observe also that the condition $E(Y|x,w)$ is independent of $w$ is a weaker condition than $Y\Perp W|X$,
 usually required for other notions of collapsibility.
For example, when $W > 0$, and
$(Y|x,w) \sim U(x - w, x + w)$, we have
$E(Y|x,w) = x. ~ \mbox{ for all}~ w. $
But,
\begin{equation}\label{eqn2.4}
F(y|x,w) = \dfrac{1}{2w} ,\hspace{4mm}     x - w  <  y  <  x + w ,
\end{equation}
showing that $Y$ and $W$ are not conditionally independent given $X$ .

Some examples for Theorem \ref{thm2.1} are the following. Suppose ($W|X=x$)$\sim N(x, 1)$ and ($Y|X=x, W=w $)$\sim N(x, w)$.
As another example, let $X>0$, ($W|X=x$)$\sim G(x, 1)$ and
($Y|X=x, W=w $)$\sim G(w, wx)$, where $G(\alpha, p)$ denote the
gamma distribution with mean  $(p/{\alpha})$. In both the cases,
$E(Y|x, w)=x$ is independent of $w$ and so   the average
collapsibility of $EDF$ $\partial E(Y|x, w)/{\partial} x$ holds.\\

\noindent We next show that condition (i) or  (ii) is only
sufficient, but not necessary. Hereafter,  $\phi(z)$ and $\Phi(z)$
 denote respectively the density and the distribution function of $Z \sim
N(0, 1).$


\begin{example} \label{exn1} {\em
Suppose ($Y|X=x, W=w$) follows uniform $U(0, (x^2+(w-x)^2))$ so
that
\begin{eqnarray}
F(y|x,w) = y(x^2+(w-x)^2)^{-1}, ~~0<y<(x^2+(w-x)^2)
\end{eqnarray}
and  $E(Y|x, w)= \frac{1}{2}(x^2+(w-x)^2).$  Assume also
$(W|X=x)\sim N(x,1)$ so that
\begin{eqnarray}
\frac{\partial}{\partial x}f(w|x) = - \phi'(w-x) = (w-x)\phi(w-x).
\end{eqnarray}
Then
\begin{eqnarray}
\int E(y|x,w)\frac{\partial}{\partial x} f(w|x)dw
&=& \frac{1}{2}\int_{-\infty}^{\infty}(x^2+(w-x)^2)(w-x)\phi(w-x)dw\nonumber\\
&=& \frac{1}{2}\bigg[x^2\int_{-\infty}^{\infty}(w-x)\phi(w-x) dw + \int_{-\infty}^{\infty} (w-x)^3\phi(w-x)dw\bigg]\nonumber\\
&=& \frac{1}{2}\bigg[x^2\int_{-\infty}^{\infty} t\phi (t) dt + \int_{-\infty}^{\infty} t^3\phi(t)dt\bigg]\nonumber\\
&=& 0, ~~\mbox{for all}~ x.
\end{eqnarray}
Thus, from (\ref{eqn2.3}), average collapsibility  over $W$ holds,
but neither condition (i) nor condition (ii) is satisfied.}
\end{example}

\noindent We next discuss an implication of Theorem \ref{thm2.1} to linear
regression models.

 \noindent {\bf Linear regression.} Consider the following conditional and marginal linear
regression models respectively:
\begin{equation*}
E(Y|X=x,W=w)= \left\{ \begin{array}{cl}
\alpha(w) + \beta(w)x, & \textrm{if $W$ is discrete}\\
\alpha + \beta x + \gamma w, &\textrm{if $W$ is continuous}
\end{array}\right.
\end{equation*}
and
\begin{equation*}
E(Y|x) = \tilde{\alpha} + \tilde{\beta}x.
\end{equation*}
Then
\begin{equation*}
\frac{\partial}{\partial x}E(Y\mid X=x,W=w)= \left\{
\begin{array}{cl}
 \beta(w), & \textrm{if $W$ is discrete}\\
 \beta , &\textrm{if $W$ is continuous}
\end{array}\right.
\end{equation*}
and
\begin{equation*}
\frac{\partial}{\partial x}E(Y\mid x) = \tilde{\beta}.
\end{equation*}

We say that the regression coefficient $\beta(w)$ (or $\beta$) is
simply collapsible if $\beta(w)$ = $\tilde{\beta} ~\mbox {for all}~
 w$ (or $\beta$ = $\tilde{\beta}$). Also, it is said to be average
collapsible if
\begin{equation}\label{eqn8n}
E_{W|x}(\beta(W)) = \tilde{\beta} ~~ (\mbox{or}~  E_{W|x}(\beta) = \tilde{\beta}), ~ \mbox{for  all}~ x.
\end{equation}
Thus, the average   collapsibility of $EDF$ reduces to the average
collapsibility of regression coefficients, in the case of linear regression models.

The average collapsibility of
regression coefficients $\beta(w)$ under the condition $E_{W}(\beta(W)) =
\tilde{\beta}$ has been discussed by  Vellaisamy and Vijay (2007).
However, the definition of average collapsibility given in
(\ref{eqn8n}) is more natural as it involves the joint
distribution of $W$ and $X$. Note also that $E_{W|x}(\beta(W)) =
\tilde{\beta}$
 for all $x$ implies $E_{W}(\beta(W)) = \tilde{\beta}$, but not necessarily conversely.

 \noindent Next, we look at the average   collapsibility of mixed
derivative of interaction ($MDI$). Since
\begin{equation}\label{neqn9}
 \frac{\partial^2}{\partial x \partial y}\log f(x,y) = \frac{\partial^2}{\partial x \partial y}\log
 f(y|x), ~~ \mbox{ for ~ all $x$ and $y$},
\end{equation}
it follows from Proposition 3.2.1 of Whittaker (1990) that
\begin{equation*}
 \frac{\partial^2}{\partial x \partial y}\log
 f(y|x)= 0   ~~ \mbox{ for ~ all $x$ and $y$}  \Longleftrightarrow Y \Perp X.
\end{equation*}

\noindent In view of \eqref{neqn9}, the $MDI$ henceforth stands
for ${\partial^2 \log
 f(y|x)}/{\partial x \partial y},$
which motivates the following definition of average
collapsibility.
 {\definition The $MDI$ is said to be average
collapsible over $W$ if
\begin{equation*}
 E_{W|x}\left(\frac{\partial^2}{\partial x \partial y}\log f(y|x,W)\right) =
 \frac{\partial^2}{\partial x \partial y}\log f(y|x),\hspace{4mm} \mbox{for all} ~  (y, x).
\end{equation*}}
 It is assumed that  $
\log f(y|x)$ has continuous partial derivatives so that
\begin{equation*}
\frac{\partial^2}{\partial x \partial y}\log f(y|x)=
\frac{\partial^2}{\partial y \partial x}\log f(y|x) ~~ \mbox{for
all} ~  (y, x). \end{equation*} The following result provides a set
of sufficient conditions for the average   collapsibility of
$MDI.$

\begin{theorem} \label{thm2.2} The $MDI$ is average collapsible over $W$ if either \\
$ \hspace*{3cm} (i)~ Y \Perp W|X$, or  \\
$ \hspace*{3cm}  (ii)~ X \Perp W|Y$\\
holds.
\end{theorem}

\noindent Xie \emph{et al.~}(2008)) showed that  condition (i)
or (ii) in Theorem  \ref{thm2.2} is necessary and sufficient for
uniform collapsibility. The following counter-example shows that
they are only sufficient, but not necessary for average
collapsibility.

\begin{example} \label{exam4n} {\em
 Let $X>0$ and $(W|x)\sim N(x, 1).$ Assume  that
\begin{equation} \label{eqnexam2n}
 f(y|x, w)=  x y^{x-1} ( x^2 + (w-x)^{2}), ~ 0<y< ( x^2 +
 (w-x)^{2})^{-1/x},
\end{equation}
which can easily be seen to be  a valid density.

Then
\begin{eqnarray} \label{eqnexam4n}
\frac{\partial^2}{\partial x\partial y}\log f(y|x, w) =
\frac{1}{y}
 = E_{W|x}\left( \frac{\partial^2}{\partial x\partial y}\log
f(y|x, W)   \right).
\end{eqnarray}

\noindent Since $(W|x)\sim N(x, 1)$, it follows that the marginal
density of $(Y|x)$ is
\begin{eqnarray*} \label{eqnexam5n}
 f(y|x)&=& \int_{-\infty}^{\infty} f(y|x, w) f(w|x) dw \\
 &=& xy^{x-1}\left[ \int_{-\infty}^{\infty} x^2 \phi(w-x)dw + \int_{-\infty}^{\infty}(w-x)^2 \phi(w-x)dw
 \right]\\
&=& x y^{x-1} (x^2 +1),
\end{eqnarray*}
which is also a valid density on $0<y< ( x^2 + 1)^{-1/x}.$

\noindent Also, it follows from \eqref{eqnexam4n}
\begin{eqnarray*} \label{eqnexam7n}
\frac{\partial^2}{\partial x\partial y}\log f(y|x)=
\frac{1}{y}=E_{W|x}\left( \frac{\partial^2}{\partial x\partial
y}\log f(y|x, W)\right).
\end{eqnarray*}
Thus, average collapsibility holds, though the condition  (i) is
not satisfied.}
\end{example}

 \noindent It was quite challenging to construct  Example \ref{exam4n}, as it requires the interchange
  of $\log$ and  integration, in addition to the other conditions. Observe also that in  Example \ref{exam4n},
\begin{eqnarray*}
 \frac{\partial}{\partial y}\log f(y|x, w)=\frac{\partial}{\partial y}\log
 f(y|x), \mbox{ for  all}~ (y,x),
\end{eqnarray*}
which leads to the average collapsibility. This observation leads to the following result which
generalizes Theorem \ref{thm2.2} whose proof is immediate.

\begin{theorem} \label{thm2.2new} The $MDI$ is average collapsible over $W$ if either
\begin{eqnarray*}
 &(i)&  \frac{\partial}{\partial y}\log f(y|x, w)=\frac{\partial}{\partial y}\log
 f(y|x), \mbox{ for  all}~ (y,x), ~~\mbox{or} \\
 &(ii)&
 \frac{\partial}{\partial x}\log f(y|x, w)=\frac{\partial}{\partial x}\log
 f(y|x), \mbox{ for  all}~ (y,x)
\end{eqnarray*}
holds.
\end{theorem}

\noindent As additional examples for Theorem \ref{thm2.2new}, let  $f(y|x, w)$ be as in Example \ref{exam4n}, consider, for $\lambda >0$, the tempered normal density
\begin{eqnarray*}
 t_{\lambda}(w|x)= c_{\lambda}(x) e^{-\lambda w} \phi(w-x)   , \mbox{ for }~ x>0, w \in \mathbb{R},
 \end{eqnarray*}
where
\begin{eqnarray*}
 c_{\lambda}(x)&=&  \Big( \int_{-\infty}^{\infty}  e^{-\lambda w} \phi(w-x) dw   \Big)^{-1}
              = e^{ ( x^2- (x-\lambda)^2)/2}.
 \end{eqnarray*}
\noindent That is,  $t_{\lambda}(w|x)= \phi(w-x+ \lambda)$. Then the corresponding marginal density of $(Y|x)$ is
\begin{eqnarray*} \label{eqnexam5nn}
 f_{\lambda}(y|x)
 &=& xy^{x-1}\left[ \int_{\infty}^{\infty} x^2 \phi(w-x+\lambda)dw + \int_{\infty}^{\infty}(w-x)^2 \phi(w-x+\lambda)dw
 \right]\\
&=& x y^{x-1} (x^2 + {\lambda}^2+1),
\end{eqnarray*}
which is also a valid density on $0<y< ( x^2 + {\lambda}^2+1)^{-1/x}.$
Thus, the average collapsibility of $MDI$ holds for the family $\{ t_{\lambda}(w|x) \}, \lambda >0,$  also.

\noindent Next, we discuss the connection of Theorem \ref{thm2.2}
to logistic regression models.

 \noindent {\bf Logistic regression}. Let $Y$ be binary and consider the following conditional and marginal logistic regression models
 (Vellaisamy and Vijay (2007), Xie \emph{et al.} (2008)) considered in the literature:
\begin{equation*}
\log\left(\frac{f(1|x,w)}{f(0|x,w)}\right)= \left\{
\begin{array}{cl}
 \alpha(w) + \beta(w)x, & \textrm{if $W$ is discrete}\\
\alpha + \beta x + \gamma w, &\textrm{if $W$ is continuous}
\end{array}\right.
\end{equation*}
and
\begin{equation*}
\log\left(\frac{f(1|x)}{f(0|x)}\right)= \tilde{\alpha} +
\tilde{\beta}x.
\end{equation*}
We say the logistic regression coefficient is simply collapsible if
\begin{equation*}
\tilde{\beta}= \left\{
\begin{array}{cl}
 \beta(w) \textrm{ for all}~ w, & \textrm{if $W$ is discrete}\\
 \beta , &\textrm{if $W$ is continuous}.
\end{array}\right.
\end{equation*}
Also, we say $\beta(w)$ or $\beta$ is said to be average
collapsible
 if
$E_{W|x}(\beta(W)) = \tilde{\beta},$ when $W$ is discrete and $E_{W|x}(\beta) = \tilde{\beta},$ when $W$ is continuous.\\
Since $Y$ is binary, the partial derivative is replaced by the difference between the adjacent levels of $Y$ (see Cox (2003)) so that
\begin{eqnarray}
\frac{\partial}{\partial x}\left(\frac{\partial}{\partial y}\log
f(y|x, w)\right)
&=& \frac{\partial}{\partial x}(\log f(1|x,w) - \log f(0|x,w)) \nonumber \\
&=& \frac{\partial}{\partial x}\log\left(\frac{f(1|x,w)}{f(0|x,w)}\right) \nonumber\\
&=& \left\{ \begin{array}{cl}
 \displaystyle \frac{\partial}{\partial x}(\alpha(w) + \beta(w)x)=\beta(w), & \textrm{if $W$ is discrete}\\
 \displaystyle \frac{\partial}{\partial x}(\alpha + \beta x + \gamma w)=\beta,
&\textrm{if $W$ is continuous},
\end{array}\right.\nonumber
\end{eqnarray}
the logistic regression coefficients corresponding to both the cases of $W$ . \\
From Theorem \ref{thm2.2}, we now conclude that $\beta(w)$ or
$\beta$ is average collapsible if $(i) Y \Perp W|X$ or $(ii) X
\Perp W|Y$ holds.

\vone \noindent Finally, we discuss a new measure called log-expectation dependence  ($LED$)  between
$X$ and $Y
>0$, defined by ${\partial}\log E(Y|x,w)/{\partial x}$, where it
is assumed that $0 < E(Y|x) < \infty$, for all $x$.
 First note that for all $x$,
\begin{eqnarray*}
\frac{\partial}{\partial x}\log E(Y|x)= 0
& \Longleftrightarrow & \frac{\partial}{\partial x} E(Y|x)=0 \\ &
\Longleftrightarrow & \int y \frac{\partial}{\partial x}
\left(dF(y|x)\right)=0\\
& \Longleftrightarrow & dF(y|x)=dF(y|x^{'}) ~~{\textrm{for ~ all
$y, x$ and $x^{'}$}} \\
& \Longleftrightarrow & Y \Perp X.
\end{eqnarray*}
Also, by Theorem 1 of Xie {\em{et~ al.}} (2008),
\begin{equation*}
\frac{\partial}{\partial x}\log E(Y|x)\geq 0 \Longrightarrow
\frac{\partial}{\partial x} E(Y|x) \geq 0 \Longrightarrow \rho(Y,
X) \geq 0,
\end{equation*}
where $\rho(Y, X)$ is the correlation coefficient between $Y$ and
$X$.

 \noindent Next, we discuss the collapsibility issues for the $LED$ measure and hence the following
 definition.

{\definition The LED is simple collapsible if
\begin{equation}\label{eqn2.8}
\frac{\partial}{\partial x}\log E(Y|x,w) =
\frac{\partial}{\partial x}\log E(Y|x), \hspace{3mm}\mbox{for all}
~  x \hspace{1mm} and \hspace{1mm} w
\end{equation}
and  average collapsible if
\begin{equation}\label{eqn2.9}
E_{W|x}\left(\frac{\partial}{\partial x}\log E(Y|x,W)\right) =
\frac{\partial}{\partial x}\log E(Y|x), \hspace{3mm}\mbox{for
all}~ x.
\end{equation}}
\begin{theorem}\label{thmled} The $LED$ is simple collapsible and hence average collapsible if
$E(Y|x,w)$ does not depend on  w.
\end{theorem}

\noindent We next discuss relevance of $LED$ in the context of
Poisson and negative binomial (NB) regression models.
{\bf Poisson regression}. Consider the Poisson regression model defined by
\begin{center}
$(Y|X = x,W = w)\sim Poi(\lambda(x,w)),$
\end{center}
where the mean
\begin{eqnarray}
E(Y|x,w) &=& \lambda(x,w)
=\left\{ \begin{array}{cl}
 e^{\alpha(w) + \beta(w)x}, & \textrm{if $W$ is discrete}\\
e^{\alpha + \beta x + \gamma w}, &\textrm{if $W$ is continuous}.
\end{array}\right.\nonumber
\end{eqnarray}
Then
\begin{eqnarray}
\frac{\partial}{\partial x} \left(\log E(Y|x,w)\right)&=&\left\{
\begin{array}{cl}
 \beta(w), & \textrm{if $W$ is discrete}\\
\beta, &\textrm{if $W$ is continuous}.
\end{array}\right.\nonumber
\end{eqnarray}
Let $(Y|x)\sim Poi(e^{\tilde{\alpha} +\tilde{\beta}x })$,
 the marginal  Poisson regression model, so that
\begin{equation*}
\log E(Y|x) = \tilde{\alpha} + \tilde{\beta}x;
~~\frac{\partial}{\partial x}\log E(Y|x) = \tilde{\beta}.
\end{equation*}

\noindent Then by Theorem \ref{thmled}, the average collapsibility
of Poisson regression coefficient $\beta(w)$
(or\hspace{2mm}$\beta$) holds, that is,
\begin{equation*}
E_{W|x}(\beta(W)) = \tilde{\beta}\hspace{4mm}(or\hspace{2mm}
E_{W|x}(\beta) = \tilde{\beta})
\end{equation*}
is true, when $ \lambda(x, w)$  does not depend on $w$ which in turn holds when for example
 $\gamma = 0$. Note that this does not in general mean that $Y
\Perp W|X.$

The following interesting example shows that  average collapsibility may hold,
even when $ E(Y|x, w)$ depends on $w$.

\begin{example} \label{nex7} {\em
Let $X>0$ and $(Y|x, w) \sim P( \lambda(x) w)$, where  $\lambda(x)
= exp{(\alpha + \beta x)}$. Then $ E(Y|x, w)=\lambda(x)w$ and
\begin{equation}\label{neqn18}
\frac{\partial}{\partial x}\log E(Y|x, w) = \beta =
E_{W|x}\left(\frac{\partial}{\partial x}\log E(Y|x, W) \right).
\end{equation}
Let now $(W|x) \sim G(x, x)$, the gamma distribution with mean
unity. Then it is known that
\begin{equation*}
 (Y|x) \sim NB \left( x, \frac{x}{x+\lambda (x)}\right),
\end{equation*}
the negative binomial (NB) distribution with
\begin{equation*}
 P(Y=y|x)= \frac{\Gamma(y+x)}{ y! ~\Gamma(x)} \left({\frac{x}{x+\lambda(x)}}\right)^{x}
 \left({\frac{\lambda(x)}{x+\lambda(x)}}\right)^{y}, ~
 y=0,1, \cdots  \cdot
\end{equation*}

\noindent Hence,
\begin{equation}\label{neqn19}
 E(Y|x)=\lambda (x); ~\frac{\partial}{\partial x}\log E(Y|x)=
 \beta.
\end{equation}
 Thus, from (\ref{neqn18}) and (\ref{neqn19}),  the average collapsibility holds. Note here the covariates $W$ and $X$ are not
independent.}
\end{example}

\noindent {\bf Negative binomial regression. }  Suppose  in Example
\ref{nex7} we assume in addition that the unobservable  $ W$ is independent of
$X$ and $ W \sim G(\theta, \theta)$. Then again
\begin{equation}{\label{eqn22n}}
   (Y|x) \sim NB \left( \theta, \frac{\theta}{\theta+\lambda (x)}\right);~~E(Y|x)= \lambda
   (x).
\end{equation}

\noindent The model \eqref{eqn22n} is the usual NB regression model. Thus, the
average collapsibility of the $LED$ function corresponds to that
of the NB regression coefficient $\beta$. It is interesting to note that
when the unobserved covariate $W$ follows the gamma distribution
with mean unity, the average collapsibility of the NB regression
coefficient holds, even when $W$ and $X$ are not independent (Example \ref{nex7}).
Note, however, in the negative binomial regression,
\begin{equation}
 Var(Y|x) = \lambda (x) \left(1+ \frac{\lambda (x)}{\theta}\right) > \lambda (x)= E(Y|x),
\end{equation}
unlike the Poisson regression case. Thus, whenever the data
exhibits over dispersion (variance exceeds mean),  the negative
binomial regression model is commonly used.


\vone \noindent {\bf 3. The Multivariate Case}

In this section, we consider an extension to the multivariate case. The case of
multivariate response $Y$ may be considered
 by treating one component at a time (Cox and Wermuth (2003) and  Xie \emph{et al.~}(2008)) and similarly the covariate
  $X$ may also be considered
 one component at a time, while keeping other components fixed.
 Therefore, we consider here only the case of multivariate random vector $W = (W_{1}, \ldots , W_{p})$. \\
A conditional measure of association, say,
$\frac{\partial}{\partial x}(E(Y|x, w)$ is simple collapsible over
$W$ if
\begin{equation*}
 \frac{\partial}{\partial x}(E(Y|x,w)) = \frac{\partial}{\partial x}(E(Y|x)),
 \hspace{3mm} \mbox{for all} ~  x \hspace{1mm} {\rm and} \hspace{1mm} w = (w_{1}, \ldots , w_{p}).
\end{equation*}
 and average collapsible if
 \begin{equation*}
  E_{W|x} \left(\frac{\partial}{\partial x}(E(Y|x,W))\right) = \frac{\partial}{\partial x}(E(Y|x)),
 \hspace{3mm} \mbox{for all} ~  x .
\end{equation*}
The definition of average
collapsibility of other measures of association remains the same,
except that $W$ is now a $p$- variate random vector.

 \noindent Let $W = (W_{1}, W_{2})$, where $W_{1}$ has $q$ components  and
$W_{2}$ has ($p-q$) components. We now have
the following result for the $EDF$ and $MDI$ and the corresponding
results for $LED$ follow easily  when $E(Y|x,w)$ is homogeneous
over $w$.
\begin{theorem} \label{thm3.1} Let $W_{1} \Perp W_{2}|X$ . Then the following results hold : \\
$~~~~~~~~~~$(a) The $EDF$ is  average collapsible over $W$ if
$(i)\hspace{1mm} Y \Perp W_{1}|(X, W_{2})$ and
 $(ii)\hspace{1mm} X \Perp W_{2}$ hold. \\
$~~~~~~~~~~$(b) The $MDI$ is average collapsible over $W$ if
$(i)\hspace{1mm} Y \Perp W_{1}|X$ and
 $(ii)\hspace{1mm} X \Perp W_{2}|Y$ hold.
\end{theorem}

 By symmetry, the average  collapsibility of $MDI$ holds when
$X$ and $Y$ are interchanged in conditions $(i)$ and $(ii)$ of
Part (b) of Theorem \ref{thm3.1}.
Also, Xie \emph{et al.~}(2008) established the uniform
collapsibility of $DDF$ and $EDF$ under an additional condition of
homogeneity of these measures. Thus, average  collapsibility holds
under less restrictive conditions and hence is applicable to a larger
class of conditional distributions that may  arise in practical applications. \\

 \noindent{\bf Acknowledgements.} This research is partially
supported by a DST  research grant No. SR/MS:706/10.

\vtwo
\renewcommand{\theequation}{A.\arabic{equation}}
 \setcounter{equation}{0}
\begin{center}
  APPENDIX
  \end{center}

\noindent {\it Proof of Theorem \ref{thm2.1}}.
Note that
\begin{eqnarray}\label{eqn2.2}
\frac{\partial}{\partial x}E(Y|x) &=& \frac{\partial}{\partial x}\int_{}^{}E(Y|x,w) f(w|x)dw\nonumber \\
&=&\int_{}^{}\frac{\partial}{\partial x}E(Y|x,w)f(w|x)dw
+ \int_{}^{}E(Y|x,w)\frac{\partial}{\partial x} f(w|x)dw\nonumber \\
&=&  E_{W|x} \left(\frac{\partial}{\partial x}E(Y|x,W)\right)+
\int_{}^{}E(Y|x,w)\frac{\partial}{\partial x} f(w|x)dw.
\end{eqnarray}
Hence, average   collapsibility holds  if and only if
\begin{equation}\label{eqn2.3}
 \int_{}^{}E(Y|x,w)\frac{\partial}{\partial x} f(w|x)dw = 0, \hspace{4mm} \mbox{for all} ~ x.
\end{equation}
Assume now condition $(i)$ holds so that
\begin{equation*}
E(Y|x,w) = h(x), \mbox{~for all} ~  x ~ \mbox{and} ~ {w},
\hspace{3mm}(say).
\end{equation*}
Then
\begin{eqnarray}
\int_{}^{}E(Y|x,w)\frac{\partial}{\partial x} f(w|x)dw&=& h(x)\int_{}^{}\frac{\partial}{\partial x} f(w|x)dw\nonumber \\
&=& h(x)\frac{\partial}{\partial x}\int_{}^{} f(w|x)dw \nonumber \\
&=& 0, \hspace{4mm}\mbox{for all} ~  x. \nonumber
\end{eqnarray}
Hence, average   collapsibility holds. \\
Assume next  condition $(ii)$ holds. Then obviously,
\begin{equation*}
\int_{}^{}E(Y|x,w)\frac{\partial}{\partial x} f(w|x)dw = 0,\hspace{4mm}  \mbox{for all} ~  x,
\end{equation*}
and so average   collapsibility holds again.

\noindent {\it Proof of Theorem \ref{thm2.2}} Since
\begin{equation}\label{eqn2.5}
  \frac{\partial^2}{\partial x \partial y}\log f(y|x) =
  \frac{\partial}{\partial x}\left(\frac{ \frac{\partial}{\partial y}f(y|x)}{f(y|x)}\right),
\end{equation}
average   collapsibility of $MDI$ holds if and only if
\begin{equation}\label{eqn2.6}
 E_{W|x}\left(\frac{\partial}{\partial x}\left( \frac{{\frac{\partial}{\partial y}f(y|x,W)}}{{f(y|x,W)}}
 \right)\right) =
  \frac{\partial}{\partial
 x}\left( \frac{{\frac{\partial}{\partial y}f(y|x)}}{{f(y|x)}}\right) ~~ \mbox{for
all}~  (y, x).
\end{equation}
Note that condition $(i)$ implies
\begin{equation*}
 f(y|x,w) = f(y|x)\hspace{4mm}  \Longrightarrow
 \frac{\partial}{\partial y}f(y|x,w) = \frac{\partial}{\partial y}f(y|x),\hspace{4mm} \mbox{for all} ~ (y, x, w).
\end{equation*}
Thus, equation (\ref{eqn2.6}) holds. \\
Observe also that
\begin{equation*}
 \frac{\partial^2}{\partial x \partial y}\log f(x,y) = \frac{\partial}{\partial y}\left(\frac{ \frac{\partial}{\partial x}(x|y)}{f(x|y)}\right),
\end{equation*}
which is the same as equation (\ref{eqn2.5}) with $x$ and $y$
interchanged. Thus, the condition $(ii)$ also implies the average
collapsibility of $MDI.$

\noindent {\it Proof of Theorem \ref{thmled}}
 Let
\begin{equation} \label{neqn15}
E(Y|x,w) = h_{1}(x)\hspace{4mm}\mbox{for all} ~  x \hspace{1mm}
{\rm and} \hspace{1mm} w.
\end{equation}
Then
\begin{eqnarray} \label{neqn16}
E(Y|x) = E_{W|x}(E(Y|x,W)) = E_{W|x}(h_{1}(x)) = h_{1}(x).
\end{eqnarray}
Thus, from \eqref {neqn15} and \eqref{neqn16},
\begin{equation*}
E(Y|x) = E(Y|x,w),\hspace{4mm}\mbox{for all} ~  x \hspace{1mm} {\rm
and} \hspace{1mm} w,
\end{equation*}
and hence simple collapsibility holds.\\
Also, since
\begin{equation*}
\frac{\partial}{\partial x}\log E(Y|x,w) =
\frac{\partial}{\partial x}\log E(Y|x)
\end{equation*}
the average  collapsibility also holds.

\noindent {\it Proof of Theorem \ref{thm3.1}} (a) Observe that
\begin{eqnarray}
E_{W|x}\left(\frac{\partial}{\partial
x}E(Y|x,W)\right)&=&\int_{w_{2}}\int_{w_{1}}\left(\frac{\partial}{\partial
x}E(Y|x,w)\right)\,dF(w_{1},w_{2}|x)\nonumber \\
&=& \int_{w_{2}}\int_{w_{1}}\left(\frac{\partial}{\partial
x}E(Y|x,w_{1},w_{2})\right)\,dF(w_{1}|x)\,dF(w_{2}|x),\hspace{4mm}(\because \hspace{2mm} W_{1}\Perp W_{2}|X)\nonumber \\
&=& \int_{w_{2}}\left(\int_{w_{1}}\frac{\partial}{\partial
x}E(Y|x,w_{2})\,dF(w_{1}|x)\right)\,dF(w_{2}|x),
\hspace{4mm}(\because \hspace{2mm} Y\Perp W_{1}|(X,W_{2}))\nonumber \\
&=& \int_{w_{2}}\left(\frac{\partial}{\partial x}E(Y|x,w_{2})\right)\,dF(w_{2}|x)\nonumber \\
&=& E_{W_{2}|x}\left(\frac{\partial}{\partial x}E(Y|x,W_{2})\right)\nonumber \\
&=& \frac{\partial}{\partial x}E(Y|x) \hspace{4mm}\mbox{for all}
~x, \nonumber
\end{eqnarray}
by condition $(ii)$ of $(a)$ and Theorem 1. \\

(b) First observe that
\begin{eqnarray}\label{eqn2.7}
\frac{{\partial^2}}{\partial x \partial y}\log f(x,y|w) &=&
\frac{\partial^2}{\partial x \partial y}
\log f(x,y,w_{1},w_{2})\nonumber \\
&=& \frac{\partial^2}{\partial x \partial y}\log f(y|x,w_{1},w_{2})\nonumber \\
&=& \frac{\partial^2}{\partial x \partial y}\log
f(w|x,w_{2}),\hspace{4mm}
\end{eqnarray}
since $Y\Perp W_{1}|X $.
By the assumption that $W_{1}\Perp W_{2}|X$  and (\ref{eqn2.7}),
\begin{eqnarray}
E_{W|x}\left(\frac{{\partial^2}}{\partial x \partial y}\log
f(x,y|W)\right) &=& \int_{w_{2}}\left( \int_{w_{1}}
\frac{\partial^2}{\partial x \partial y}\log f(y|x,w_{2})\,dF(w_{1}|x)\right)\,dF(w_{2}|x)\nonumber \\
&=& \int_{w_{2}}\left(\frac{\partial^2}{\partial x \partial y}\log f(y|x,w_{2})\right)\,dF(w_{2}|x)\nonumber \\
&=& E_{W_{2}|x}\left(\frac{\partial^2}{\partial x \partial y}\log f(y|x,W_{2})\right)\nonumber \\
&=& \frac{\partial^2}{\partial x \partial y}\log f(y|x)
\hspace{4mm}\mbox{for all} ~  x \hspace{1mm} {\rm and}
\hspace{1mm} y\nonumber
\end{eqnarray}
by condition $(ii)$ of Theorem 2.

\vone
 \noindent {\bf References}

\noindent
Cochran, W. G. (1938). The omission or addition of an independent
variable in multiple linear regression.
{\it J. R. Statist. Soc. Suppl.}, {\bf 5}, 171-176.\\
Cox, D. R. (2003). Conditional and marginal association for binary
random variables.
{\it Biometrika}, {\bf 90}, 982-984.\\
Cox, D. R. and Wermuth, N. (2003). A general condition for
avoiding effect reversal after marginalization.
{\it J. R. Statist. Soc. B}, {\bf 65}, 937-941.\\
Geng, Z. (1992). Collapsibility of relative risk in contingency
tables with a response variable.
{\it J. R. Statist. Soc. B}, {\bf 54}, 585-593.\\
Ma, Z., Xie, X. and Geng, Z. (2006). Collapsibility of
distribution dependence.
{\it J. R. Statist. Soc. B}, {\bf 68}, 127-133.\\
Mantel, N. and Haenszel, W. (1959). Statistical aspects of the
analysis of data from retrospective studies of disease.
{\it J. Natn. Cancer Inst.}, {\bf 22}, 719–748.\\
Simpson, E.H. (1951). The interpretation of interaction in
contingency tables. {\it J. R. Statist. Soc. B}, {\bf 13},
238-241.\\
Vellaisamy, P. and Vijay, V. (2008). Collapsibility of regression
coefficients and its extensions.
{\it J. Statist. Plann. Inference}, {\bf 138}, 982-994.\\
Vellaisamy, P. (2011). Average  collapsibility of distribution
dependence and quantile regression coefficients. To appear in {\it
Scand. J.
Statist.}\\
Wermuth, N. (1989). Moderating effects of subgroups in linear
models.
{\it Biometrika}, {\bf 76}, 81-92.\\
Wermuth, N. (1987). Parametric collapsibility and the lack of
moderating effects in contingency tables
with a dichotomous response variable. {\it J. R. Statist. Soc. B}, {\bf 49}, 353-364.\\
Whittaker, J. (1990). {\it Graphical Models in Applied
Multivariate Statistics}. Wiley, New York.\\
Whittemore, A. S. (1978). Collapsibility of multidimensional
contingency tables.
{\it J. R. Statist. Soc. B}, {\bf{40}}, 328-340.\\
Xie, X., Ma, Z. and Geng, Z. (2008). Some association measures and
their collapsibility. {\it Statist. Sinica}, {\bf 18},
1165-1183.\\
Yule, G.U. (1903). Notes on the theory of association of
attributes. {\it Biometrika}, {\bf 2}, 121-134.


\end{document}